\documentclass[11pt]{article}

\usepackage{amsmath}
\usepackage{amssymb}
\usepackage{eucal}

\usepackage{stmaryrd}

\begin{document}
                                    \baselineskip 16pt

\title{A contribution to the theory of $\sigma$-properties of a finite group\thanks{Research was supported by the National
 Natural Science Foundation of China
 (No. 12171126, 12101165).
 Research of the third author and the  fourth author was
 supported by the Ministry of Education of the
 Republic of Belarus (No.~20211328, 20211778).  
}\\{\small  In memory of Professor Francesco de Giovanni}}

\author{A-Ming Liu, \ \ Wenbin Guo  \\
{\small School of Mathematics and Statistics, Hainan University, Haikou, Hainan,  570228, P.R. China}\\
{\small  E-mail: amliu@hainanu.edu.cn, \ wbguo@ustc.edu.cn } \\ \\
Vasily G. Safonov \\ 
{\small Institute of Mathematics of the National Academy of Sciences of Belarus,}\\
{\small Minsk 220072, Belarus.} {\small E-mail: vgsafonov@im.bas-net.by}\\ 
{\small Department of Mechanics and Mathematics, Belarusian State University,} \\ 
{\small Minsk 220030, Belarus.} {\small E-mail: vgsafonov@bsu.by}\\ \\
Alexander N. Skiba\\
{\small Department of Mathematics and Technologies of Programming,}\\
{\small  Francisk Skorina Gomel State University,}\\
{\small Gomel 246019, Belarus}\\
{\small E-mail: alexander.skiba49@gmail.com}}

\date{}
\maketitle

\begin{abstract}      Let   $\sigma =\{\sigma_{i} \mid i\in I\}$ be some
 partition of the set of all primes.  
 A  subgroup $A$ of a finite group $G$ is said to be: (i)  
 \emph{$\sigma$-subnormal} in $G$ if   
   there is a subgroup chain  $A=A_{0} \leq A_{1} \leq \cdots \leq
A_{n}=G$  such that  either $A_{i-1} \trianglelefteq A_{i}$ or 
$A_{i}/(A_{i-1})_{A_{i}}$ is a  ${\sigma}_{j}$-group, $j=j(i)$,    for all $i=1, \ldots , n$; 
(ii)  \emph{modular} in $G$ 
 if the following conditions are held:
(1) $\langle X, A \cap Z \rangle=\langle X, A \rangle \cap Z$ for all $X \leq G, Z \leq
 G$ such that
$X \leq Z$, and
(2) $\langle A, Y \cap Z \rangle=\langle A, Y \rangle \cap Z$ for all $Y \leq G, Z \leq
 G$ such that  $A \leq Z$; (iii) \emph{$\sigma$-quasinormal in $G$} if $A$ is 
$\sigma$-subnormal and modular in $G$.

   We obtain a  description of  finite  groups in which 
 $\sigma$-quasinormality (respectively, modularity) is a
 transitive relation.  Some known results are extended.



\end{abstract}

\section{Introduction}

Throughout this paper, all groups are finite and $G$ always denotes
a finite group; ${\cal L}(G)$ is the lattice of all subgroups of $G$;
 $G$ is said to be an  \emph{$M$-group} \cite{Schm}
 if the lattice ${\cal L}(G)$  is   modular.   
  Moreover,  $\mathbb{P}$ is the set of all  primes, $\pi \subseteq \mathbb{P}$, 
 $\pi'=\mathbb{P}\setminus \pi$, and  $\sigma =\{\sigma_{i} \mid
 i\in I \}$ is some  partition of $\mathbb{P}$.
 If
 $n$ is an integer, the symbol $\pi (n)$ denotes
 the set of all primes dividing $n$; as usual,  $\pi (G)=\pi (|G|)$, the set of all
  primes dividing the order of $G$; $\sigma (n)= \{ \sigma_{i}\mid \sigma_{i} \cap \pi
 (n)\ne \emptyset\}$ and 
 $\sigma (G)= \sigma (|G|)$ \cite{commun, 1}.
 A group $G$ is said to be \cite{commun, 1}: \emph{$\sigma$-primary} if $G$ is a $\sigma _{i}$-group
 for some $i$;   \emph{$\sigma$-nilpotent} if $G$ is a direct product
 of $\sigma$-primary groups.

A subgroup $A$ of $G$ is said to be \emph{quasinormal} (Ore)
 or \emph{permutable} (Stonehewer)
 in $G$
 if  $A$ permutes with every subgroup $H$ of $G$, that is, $AH=HA$.

The quasinormal subgroups have many interesting and
 useful for applications  properties.
For instance, if $A$ is quasinormal in $G$, then:  {\sl $A$ is subnormal in $G$}
 (Ore \cite{5}),  {\sl $A/A_{G}$ is nilpotent} (Ito and Szep \cite{It}), 
 {\sl every chief factor 
$H/K$ of $G$ between $A_{G}$ and $A^{G}$ is central, that is,  $C_{G}(H/K)=G$} ({Maier and Schmid  \cite{MaierS}),
  and,
 in general,  the section \emph{$A/A_{G}$ is not necessarily abelian}
(Thomson  \cite{Th}).

Quasinormal subgroups have a
 close connection with the so-called modular subgroups.

Recall that a subgroup $M$ of $G$  is said to be: (i)  \emph{ modular}  in $G$ \cite{Schm}
   if $M$ is a modular element  (in the \emph{sense of
 Kurosh} \cite[p. 43]{Schm})  of the  lattice ${\cal L}(G)$,  that is,  
(1) $\langle X,M \cap Z \rangle=\langle X, M \rangle \cap Z$ for all $X \leq G, Z \leq
 G$ such that $X \leq Z$, and  
(2) $\langle M, Y \cap Z \rangle=\langle M, Y \rangle \cap Z$ for all $Y \leq G, Z \leq
 G$ such that  $M \leq Z$; (ii) \emph{submodular} in $G$ if   
   there is a subgroup chain  $A=A_{0} \leq A_{1} \leq \cdots \leq
A_{n}=G$  such that   $A_{i-1}$ is modular in $ A_{i}$   for all $i=1, \ldots , n$.

Every quasinormal is clearly modular in the group. Moreover,
 the following remarkable result is well-known.

{\bf Theorem A} (Schmidt \cite[Theorem 5.1.1]{Schm}) {\sl A subgroup $A$ of $G$ is
 quasinormal in $G$ if and only if $A$ is   subnormal and modular in $G$}.

 This  result made
 it possible to find an analogue of quasinormality in the theory of the
 \emph{$\sigma$-properties} of a group \cite{ProblemI}. 

A subgroup $A$ of $G$ is said to
 be \emph{$\sigma$-subnormal}
  in $G$ \cite{commun, 1} if  there is a subgroup chain
 $A=A_{0} \leq A_{1} \leq \cdots \leq
A_{n}=G$  such that  either $A_{i-1} \trianglelefteq A_{i}$ or
$A_{i}/(A_{i-1})_{A_{i}}$ is  ${\sigma}$-primary 
  for all $i=1, \ldots , n$; \emph{$\sigma$-seminormal
in $G$} (J.C. Beidleman) if $x\in N_{G}(A)$ for all $x\in G$ such that
 $\sigma (|x|)\cap \sigma (A)=\emptyset$.

{\bf Definition 1.1.} We say that a subgroup $A$ of $G$
 is \emph{$\sigma$-quasinormal} in $G$ if  $A$ is  $\sigma$-subnormal and modular
 in $G$.

Before continuing, consider some examples.

{\bf Example 1.2.}  (i) 
In the first limiting case, when 
  $\sigma =\{\mathbb{P}\}$, every group  is $\sigma$-nilpotent and every 
subgroup of any group is $\sigma$-subnormal. Therefore in this case
 a subgroup $A$ of $G$ is  $\sigma$-quasinormal  if and only if it is modular in $G$.

(ii)    In the second limiting case, when   $\sigma =
\sigma ^{1}=\{\{2\}, \{3\}, \{5\} \ldots  \}$,  a subgroup $A$ of $G$
 is  $\sigma$-subnormal in $G$   if and only if it is subnormal in $G$.
Therefore in this case, in view of Theorem A,  
 a subgroup $A$ of $G$ is  $\sigma$-quasinormal  if and only if it is quasinormal 
 in $G$. 

(iii)  In the case
 $\sigma =\sigma ^{1\pi }=\{\{p_{1}\}, \ldots , \{p_{n}\}, \pi'\}$, where  
$\pi= \{p_{1}, \ldots , p_{n}\}$, a subgroup $A$ of
 $G$ is $\sigma^{1\pi }$-subnormal in $G$ if and only if  $G$ has 
 a subgroup chain
 $A=A_{0} \leq A_{1} \leq \cdots \leq
A_{n}=G$  such that  either $A_{i-1} \trianglelefteq A_{i}$  or
$A_{i}/(A_{i-1})_{A_{i}}$ is a  ${\pi}'$-group   for all $i=1, \ldots , n$.

In this case we say, following 
 \cite{comm, ???, ????}, that $A$ is \emph{$1 \pi$-subnormal} in $G$, and we say that 
$A$ is \emph{$1 \pi$-quasinormal} in $G$ if $A$ is 
$1 \pi$-subnormal and modular in $G$.   Note, in passing,  that  $A$ is   $1 \pi$-subnormal
 in $G$ if and only if $A$ is \emph{$\mathfrak{F}$-subnormal in $G$ in
 the sence of Kegel}
 \cite{KegI}, where $\mathfrak{F}$ is the class of all $\pi'$-groups.

(iv) In the other classical  case $\sigma =\sigma ^{\pi}=\{\pi,
\pi'\}$  a
 subgroup $A$ of  $G$  is
 ${\sigma} ^{\pi}$-subnormal  in $G$ if and only
 if  $G$ has   a subgroup chain
 $A=A_{0} \leq A_{1} \leq \cdots \leq
A_{n}=G$  such that  either $A_{i-1} \trianglelefteq A_{i}$, or
$A_{i}/(A_{i-1})_{A_{i}}$ is a ${\pi}$-group, or  
 $A_{i}/(A_{i-1})_{A_{i}}$ is  a ${\pi}'$-group  for all $i=1, \ldots , n$.

In this case we  say that $A$ is
 \emph{$\pi, \pi'$-subnormal} in $G$ \cite{comm, ???, ????}, and we say that 
$A$ is \emph{$\pi, \pi'$-quasinormal} in $G$ if $A$ is 
$\pi, \pi'$-subnormal and modular in $G$.

The
 theory of $\sigma$-quasinormal subgroups was constructed in the paper \cite{Hu11}.
 In particular,  it was proven  the following result covering in the
 case $\sigma =
\sigma ^{1}=\{\{2\}, \{3\}, \{5\} \ldots  \}$ the above mentioned results in
 \cite{5, It, MaierS}.

{\bf Theorem B} (See Theorem C in \cite{Hu11}). {\sl 
 Let $A$ be a $\sigma$-quasinormal subgroup of   $G$.
 Then the following statements hold:}

(i) {\sl  $A$ permutes with all  Hall $\sigma _{i}$-subgroups  of  $G$ for all $i$. }

(ii) {\sl The quotients
$A^{G}/A_{G}$ and $G/C_{G}(A^{G}/A_{G}) $ are $\sigma$-nilpotent, and }

(iii) {\sl Every chief factor $H/K$  of $G$ between
 $A^{G}$ and $A_{G}$ is $\sigma$-central in $G$ , that is,
 $(H/K)\rtimes  (G/C_{G}(H/K)) $ is $\sigma$-primary.  }

(iv)  {\sl For every $i$ such that $\sigma _{i} \in \sigma
(G/C_{G}(A^{G}/A_{G}))$  we have
 $\sigma _{i} \in  \sigma (A^{G}/A_{G}).$
}

(v) {\sl $A$  is $\sigma$-seminormal in $G$.}

A group $G$ is said to be a
\emph{$PT$-group} \cite[2.0.2]{prod}
 if  quasinormality is
a transitive relation on $G$, that is,    if $H$ is a quasinormal subgroup 
 of $K$ and $K$ is a quasinormal subgroup 
 of $G$, then $H$ is a quasinormal subgroup of $G$.

The   description of  $PT$-groups   was first obtained by    Zacher \cite{zaher},
  for the soluble  case, and   by Robinson in \cite{217}, for the general case.

Bearing in mind  the results in
 \cite{zaher, 217} and many other known results on $PT$-groups (see, in particular, 
 Chapter 2 in \cite{prod}), it seems to be  natural to ask:

{\bf Question 1.3.} {\sl What is the structure of  $G$
 provided   $\sigma$-quasinormality is  a transitive relation in $G$?}

{\bf Question 1.4.} {\sl What is the structure of  $G$
 provided    modularity  is a transitive relation in $G$?}

Note that in view of Example 1.2(i), Question 1.4
 is a special case of Question 1.3, where  $\sigma =\{\mathbb{P}\}$. Note also that
 for the case when $G$  is a soluble group,
 the answers to both of these questions are known.

 Frigerio proved \cite{A. Frigerio} (see also \cite{mod}) that modularity 
 is a transitive relation in a soluble group $G$ if and only if $G$ is an $M$-group.

An important step in solving the general Problem 1.3
 was made in the paper \cite{MZ}, where it was proven the following theorem  
 turn into Frigerion result  in the case where $\sigma =\{\mathbb{P}\}$.

{\bf Theorem C  } (X.-F. Zhang,  W. Guo,   I.N. Safonova,  
A.N. Skiba  \cite{MZ}).  {\sl       Let   $G$ be  a  soluble 
 group   and  $D=G^{\frak{N_{\sigma}}}$. Then $\sigma$-quasinormality
 is  a transitive relation in $G$ if and only if 
    the following conditions hold:} 

(i) {\sl $G=D\rtimes M$, where $D$   is an abelian  Hall
 subgroup of $G$ of odd order, $M$ is a $\sigma$-nilpotent $M$-group.}

(ii) {\sl  every element of $G$ induces a
 power automorphism on  $D$,   }

(iii) {\sl  $O_{\sigma _{i}}(D)$ has 
a normal complement in a Hall $\sigma _{i}$-subgroup of $G$ for all $i$.}

{\sl Conversely, if  Conditions (i), (ii)
 and (iii) hold for  some subgroups $D$ and $M$ of
 $G$, then  $\sigma$-quasinormality
 is  a transitive relation in $G$.}

In this theorem,     $G^{\frak{N_{\sigma}}}$  denotes the \emph{$\sigma$-nilpotent
 residual} of $G$,
 that is,  the intersection of all normal subgroups $N$ of $G$ with 
$\sigma$-nilpotent quotient $G/N$.

{\bf Definition 1.5.} We say that $G$ is:  (i) a \emph{$Q\sigma T$-group} if the 
 $\sigma$-quasinormality is
a transitive relation on $G$, that is,    if $H$ is a $\sigma$-quasinormal subgroup 
 of $K$ and $K$ is a $\sigma$-quasinormal subgroup 
 of $G$, then $H$ is a $\sigma$-quasinormal subgroup of $G$; 
 (ii)  an \emph{$M T$-group}
 if  the modularity  is a transitive relation in $G$.
 
It is clear that an $M T$-group is exactly a $Q\sigma T$-group 
 where $\sigma =\{\mathbb{P}\}$.

In this article,  expanding the corresponding results of the papers 
 \cite{217, MZ, Archiv},  we answer Questions 1.3 and 1.4 in the general case.

 {\bf Definition 1.6.}   We say    that
$(D, Z(D); U_{1},  \ldots , U_{k})$  is  a \emph{Robinson complex}
 if the following fold:

(i)  $D\ne 1$ is a perfect normal subgroup of $G$,

(ii) $D/Z(D)=U_{1}/Z(D)\times \cdots \times U_{k}/Z(D)$, where $U_{i}/Z(D)$ is a  
 simple  non-abelian chief factor of $G$, $Z(D)=\Phi (D)$,  and

(iii) every chief factor of $G$  below $Z(D)$ is cyclic.

{\bf Example 1.7.}  Let  $G=SL(2, 7)\times A_{7}\times A_{5}\times B$,
 where  $B=C_{43}\rtimes C_{7}$
is  a non-abelian group of order 301.   
  Then 
 $$(SL(2, 7)\times A_{5}\times A_{7}, Z(SL(2, 7)); SL(2, 7), A_{5}Z(SL(2, 7)),
 A_{7}Z(SL(2, 7))))$$  is a  Robinson complex of $G$.

Now let $G=A_{n}\wr C_{p}=K\rtimes C_{p}$, where $K$ is the base group of the regular
 wreath product of the alternating  group $A_{n}$  of degree $n > 4$
 with a group $C_{p}$  of prime order $p$.  Then  $K$ is a minimal normal subgroup of $G$
 by \cite[Chapter A, 18.5(a) ]{DH}. Hence $G$ has no a  Robinson complex.

We say, following Robinson \cite{217},  that \emph{$G$  satisfies}:

(1)  ${\bf N}_{p}$  if whenever $N$ is  a soluble normal
subgroup of $G$, $p'$-elements of $G$ induce power automorphism in
  $O_{p}(G/N)$;

(2)  ${\bf P}_{p}$  if whenever $N$ is  a soluble normal
subgroup of $G$, every subgroup of $O_{p}(G/N)$ is quasinormal
 in  every  Sylow  $p$-subgroup of $G/N$.

 Every subnormal subgroup is both submodular and $\sigma$-subnormal in the group.
 Thus  the following well-known 
 result partially  describes the structure of  insoluble $Q\sigma T$-groups.

{\bf Theorem D}  (Robinson  \cite{217}).  {\sl $G$ is a $PT$-group if 
 and   only if  $G$  has a normal perfect subgroup $D$ such that:}

(i) {\sl  $G/D$ is a soluble $PT$-group, and }

(i) {\sl if $D\ne 1$, $G$ has a Robinson complex
 $(D, Z(D); U_{1},  \ldots , U_{k})$ and }

(iii) {\sl   for any set  $\{i_{1}, \ldots , i_{r}\}\subseteq \{1, \ldots , k\}$, where
 $1\leq r  < k$,  $G$ and $G /U_{i_{1}}'\cdots U_{i_{r}}'$ satisfy
 ${\bf N}_{p}$ for all $p\in \pi (Z(D))$ and
 ${\bf P}_{p}$ 
for all $p \in \pi (D)$. }

Now, recall that $G$ is a non-abelian $P$-group (see \cite[p. 49]{Schm}) if 
 $G=A\rtimes \langle t \rangle$, where $A$ is    an elementary abelian
$p$-group and an element $t$ of
 prime order $q\ne p$ induces a non-trivial power
 automorphism  on $A$. In this case we say that $G$ is a \emph{$P$-group
 of type  $(p, q)$}.

{\bf Definition 1.8.}  We say that:

(i)   $G$  \emph{satisfies} 
 ${\bf Q }_{\sigma (p, q)}$    if  
 whenever $N$ is  a soluble normal
subgroup of $G$ and $P/N$ is a normal $\sigma$-primary 
 $P$-subgroup  of type  $(p, q)$  of $G/N$, 
 every  subgroup of $P/N$ is   modular in   $G/N$.

If $G$  satisfies 
 ${\bf Q }_{\sigma (p, q)}$ and  $\sigma =\{\mathbb{P}\}$,
 then say, following \cite{Archiv}, that   $G$  satisfies 
${\bf M }_{p, q}$.

(ii)   $G$  \emph{satisfies} 
 ${\bf Q }_{\sigma P}$    if  $G$  satisfies 
 ${\bf Q }_{\sigma (p, q)}$ for each  pair $p, q$
 such that there is a  $P$-group  of type  $(p, q)$.

    In this paper, based on Theorems C and   D, 
 we prove the following result.

{\bf Theorem E.}  {\sl A group $G$ is a $Q\sigma T$-group if   
 and   only if  $G$  has a perfect normal subgroup $D$ such that:}

(i) {\sl  $G/D$ is a soluble  $Q\sigma T$-group,  }

(ii) {\sl if  $D\ne 1$,  $G$ has a Robinson complex
 $(D, Z(D); U_{1},  \ldots , U_{k})$ and }

(iii) {\sl   for any set $\{i_{1}, \ldots , i_{r}\}\subseteq \{1, \ldots , k\}$, where
 $1\leq r  < k$,  the groups $G$ and $G /U_{i_{1}}'\cdots U_{i_{r}}'$ satisfy
 ${\bf N}_{p}$ for all $p\in \{2, 3\} \cap \pi (Z(D))$,
 ${\bf P}_{p}$ 
for all $p\in \pi (D)$, and ${\bf Q }_{\sigma (p, q})$
  for all 
 $\{p, q\}\cap  \pi (D)\ne \emptyset$.}

 Theorem E gives a solution to Question 1.3. The following special case  
of Theorem E gives a solution to Question 1.4. 
                                      
{\bf Theorem F. }   {\sl
 A group $G$ is an $MT$-group if   
 and   only if  $G$  has a perfect normal subgroup $D$ such that:}

(i) {\sl  $G/D$ is an $M$-group,  }

(ii) {\sl if  $D\ne 1$,  $G$ has a Robinson complex
 $(D, Z(D); U_{1},  \ldots , U_{k})$ and }

(iii) {\sl   for any set $\{i_{1}, \ldots , i_{r}\}\subseteq \{1, \ldots , k\}$, where
 $1\leq r  < k$,  $G$ and $G /U_{i_{1}}'\cdots U_{i_{r}}'$ satisfy
 ${\bf N}_{p}$ for all $p\in \{2, 3\}\cap \pi (Z(D))$,
 ${\bf P}_{p}$ 
for all $p\in \pi (D)$, and ${\bf M}_{p, q}$ for all pairs
 $\{p, q\}\cap \pi (D)\ne \emptyset.$}

We prove Theorem E (and so Theorem F, as well) in Section 3. In Section 4
 we discuss some other applications of these results.

\section{Preliminaries}

The first lemma is a corollary of  general properties  of modular
 subgroups \cite[p. 201]{Schm} and $\sigma$-subnormal subgroups \cite[Lemma 2.6]{1}.

{\bf Lemma 2.1.}  {\sl Let $A$, $B$  and
$N$ be subgroups of $G$, where $A$ is $\sigma$-quasinormal and $N$ is
normal in $G$.}

(1) {\sl The subgroup $A\cap B$ is  $\sigma$-quasinormal  in   $B$.}

(2) {\sl The subgroup  $AN/N$
 is $\sigma$-quasinormal  in  $G/N$}.

(3) {\sl If $N\leq B$ and $B/N$ is  $\sigma$-quasinormal  in  $G/N$, then 
$B$  is
$\sigma$-quasinormal  in  $G$. } 

(4)   {\sl  $B$ is $\sigma$-quasinormal  in $G$, then $\langle A, B \rangle$
is $\sigma$-quasinormal  in $G$.}

{\bf Lemma 2.2} {\sl A subgroup $A$ of $G$ is a 
 maximal $\sigma$-quasinormal subgroup of
 $G$
 if and only if either  $A$ is normal in $G$ and   $G/A$ is a  
 simple gropup   or $A_{G} < A$ and 
 $G/A_{G}$ is a $\sigma$-primary non-abelian group  of order $pq$ for primes
 $p$ and $q$.}

{\bf Proof. }  First assume that  $A$ is a maximal $\sigma$-quasinormal subgroup of 
 $G$. If $A$ is  normal in $G$, then $G/A=G/A_{G}$ is  
 simple. Now assume that $A$ is not normal in $G$, so $A^{G}=G$ and, in view of
 Theorem B(ii),
 $G/A_{G}$ is a $\sigma _{i}$-group for some $i$. Hence  every subgroup of $G$ 
containing  $A_{G}$  is $\sigma$-subnormal in $G$ by \cite[Lemma 2.6(5)]{1}.
On the other hand, $U/A_{G}$ is modular in $G$ if and only if $U$ is modular in $G$
 by \cite[Page 201, Properties  (3)(4)]{Schm}. Therefore, in fact, $A$ is a
 maximal modular subgroup of $G$. 
 Hence    $G/A_{G}$ is a non-abelian group   of order $pq$ for
 primes $p, q \in \sigma _{i}$  by \cite[Lemma 5.1.2]{Schm}.

Now assume that $A_{G} < A < G$ and $G/A_{G}$ is a $\sigma$-primary non-abelian group 
  of order $pq$ for primes
 $p$ and $q$. Then $A$ is a maximal subgroup of $G$ and $A$ is a 
 $\sigma$-subnormal subgroup of $G$. Moreover,   $A/A_{G}$   is  modular
 in $G/A_{G}$ by \cite[Lemma 5.1.2]{Schm},
 so $A$ is a maximal modular subgroup of $G$
 by \cite[Page 201, Property (4)]{Schm}. Hence  $A$  is a maximal $\sigma$-quasinormal subgroup of 
 $G$.

           Finally, assume that  $A$ is normal in $G$ and 
 $G/A$ is a 
 simple non-abelian  group, then $A$ is a maximal modular
 subgroup of $G$ by \cite[Lemma 5.1.2]{Schm} and $A$
  $\sigma$-subnormal in $G$. Hence  $A$  is a maximal $\sigma$-quasinormal subgroup of 
 $G$.  
The lemma is proved.

We say that a  subgroup $A$ of $G$ is said to be  \emph{$\sigma$-subquasinormal} in $G$ 
if     
   there is a subgroup chain  $A=A_{0} \leq A_{1} \leq \cdots \leq
A_{n}=G$  such that  $A_{i-1}$ is $\sigma$-quasinormal in $ A_{i}$ 
 for all $i=1, \ldots , n$.

It is clear that $G$ is a $Q\sigma T$-group if and only if every of its 
$\sigma$-subquasinormal subgroups is $\sigma$-quasinormal in $G$.

The class of groups $\frak{F}$  is a \emph{hereditary formation} if 
$\frak{F}$   is closed under
 taking  derect products, homomorphic images and  subgroups.
If $\frak{F}\ne \emptyset$  is a hereditary formation,
 then the symbol $G^{\frak{F}}$   denotes the \emph{$\frak{F}$-residual
 of $G$}, that is, the intersection of all normal subgroups $N$ of $G$ with
  $G/N\in \frak{F}$.

We use $\frak{A}^{*}$ to denote the class of all abelian groups of squarefree exponent.
It is clear that $\frak{A}^{*}$ is a hereditary formation.

{\bf Lemma 2.3.} {\sl Let  $A$,  $B$ and $N$ be subgroups of $G$, where 
  $A$
is $\sigma$-subquasinormal $G$ and $N$ is normal $G$  in $G$.  }

(1) {\sl $A\cap B$    is  $\sigma$-subquasinormal $G$  in   $B$}.

(2) {\sl $AN/N$ is  $\sigma$-subquasinormal $G$  in $G/N$. }

(3) {\sl If $N\leq K$ and $K/N$ is $\sigma$-subquasinormal $G$ 
in $G/N$, then $K$ is  $\sigma$-subquasinormal $G$  in $G.$}

 (4)    {\sl  $A^{{\frak{A}^{*}}}$    is subnormal in $G$. }

 (5)    {\sl  If $G=U_{1}\times \cdots \times U_{k}$, where $U_{i}$ is a simple
 non-abelian group, then $A$ is normal in $G$. }

{\bf Proof.}  (1)--(4). These assertions follow from  Lemma 2.6 in \cite{1} and 
corresponding lemmas in \cite{mod}.
                                   
(5)  Let $E=U_{i}A$, where $U_{i}\nleq A$.  We show that $A  \trianglelefteq  E$.

The subgroup  $A$ is $\sigma$-subquasinormal $G$ in $E$ by Part (1) and $A < E$, so
 there is a subgroup chain 
$A=E_{0} < E_{1} < \cdots < E_{t-1} < E_{t}=E$ such that $E_{i-1}$  is a maximal
$\sigma$-quasinormal subgroup of  $E_{i}$ for all $i=1, \ldots, t$ and for $M=E_{t-1}$
 we have $M=A(M\cap U_{i})$, where $M\cap U_{i}$ is 
 $\sigma$-subquasinormal in $U_{i}$. Then $M\cap U_{i}   < U_{i}$
 since $M < E$. Therefore $M\cap U_{i}=1=A\cap U_{i}$ by Lemma 2.2 since $U_{i}$ is a simple
 non-abelian group, so $A$ is a maximal $\sigma$-quasinormal subgroup of  $E$.   
Assume that $A$ is not normal in $E$. Then $E/A_{E}=U_{i}A/A_{E}$ is a 
 group of  order $qr$ for primes  $q$ and $r$ by Lemma 2.2,  where 
$U_{i}\simeq U_{i}A_{E}/A_{E}\leq E/A_{E}$. This contradiction show that
 $U_{i}\leq N_{E}(A)$, so $G\leq N_{G}(A)$. Hence we have (5).   
The lemma is proved.

{\bf Lemma 2.4.}   {\sl  If  $G$ is a 
 $Q\sigma T$-group, then every   quotient $G/N$ of $G$ is also a
$ Q\sigma T$-group. }

{\bf Proof.}  
Let $L/N$ be a $\sigma$-subquasinormal subgroup of $G/N$. Then $L$ is a
$\sigma$-subquasinormal subgroup in $G$ by Lemma 2.3(3), so $L$ 
 is $\sigma$-quasinormal in $G$ by hypothesis 
 and
 then $L/N$ is $\sigma$-quasinormal in $G/N$ by 
  Lemma 2.1(2). Hence 
$G/N$  is  a $Q\sigma T$-group.  
 The lemma is proved.

{\bf Lemma 2.5.}   {\sl If $G$ is a $Q\sigma T$-group, then $G/R$   satisfies
 ${\bf Q }_{\sigma P}$  for  every normal subgroup $R$ of $G$.}

{\bf Proof.}  
 In view of Lemma 2.4, we can assume without loss of
 generality that $R=1$.

 Let $P/N$ be any normal $\sigma$-primary non-abelian  $P$-subgroup  
 of type  $(p, q)$ of $G/N$    
  and let $L/N\leq P/N$.
 Then $L/N$ is modular in $P/N$ by \cite[Lemma 2.4.1]{Schm}, so 
 $L/N$ is submodular in $G/N$. On the other hand, $L/N$ is $\sigma$-subnormal in $G/N$ 
 since $P/N\leq O_{\sigma _{i}}(G/N) $ for some $i$.
Therefore $L/N$ is $\sigma$-subquasinormal in $G/N$ and so  
 $L$ is $\sigma$-subquasinormal in $G$  by Lemma 2.3(3). 
Hence  $L$  is 
 $\sigma$-quasinormal in $G$ by hypothesis, so $L/N$ is modular in $G/N$ by
  \cite[Page 201, Property (3)]{Schm}.   
Therefore $G$   satisfies ${\bf Q }_{\sigma P}$.    The
 lemma is proved.

We use $G^{\mathfrak{S}}$  (respectively, $G^{\mathfrak{U}}$) to denote the soluble
 (respectively, the supersoluble)  residual of $G$.

{\bf Lemma 2.6.}  {\sl Let $G$ be a  non-soluble group  
 and suppose that $G$ has a 
  Robinson complex
 $(D, Z(D); U_{1}, \ldots ,
U_{k}),$ where $D=G^{\mathfrak{S}}=G^{\mathfrak{U}}$.
  Let $U$ be a  $\sigma$-subquasinormal non-$\sigma$-quasinormal
  subgroup of $G$ of minimal order.   Then:}

(1) {\sl If $UU_{i}'/U_{i}'$ is $\sigma$-quasinormal
 in  $G/U_{i}'$ for
all $i=1, \ldots, k$, then $U$ is supersoluble.}

(2) {\sl If $U$ is  supersoluble and $UL/L$ is $\sigma$-quasinormal
 in  $G/L$ for
  all non-trivial nilpotent  normal subgroups $L$ of $G$, then
$U$ is a cyclic $p$-group for some prime $p$. }

{\bf Proof. }  Suppose that this lemma is false and let $G$ be a 
counterexample of minimal order.

(1)  Assume this is false. Suppose that 
$U\cap D\leq Z(D)$. Then every chief factor of $U$ below
 $U\cap Z(D)=U\cap D$ is cyclic  and, also,  $UD/D\simeq U/(U\cap  D)$ is
 supersoluble.
 Hence  $U$ is supersoluble, a contradiction. Therefore
 $U\cap D\nleq Z(D)$.   Moreover, Lemma 2.3(1)(2) 
 implies that $(U\cap D)Z(D)/Z(D)$ is
$\sigma$-subquasinormal in $D/Z(D)$ and so  $(U\cap
D)Z(D)/Z(D)$ is a non-trivial
normal  subgroup of  $D/Z(D)$  by Lemma 2.3(5).

  Hence for some $i$ we 
have $U_{i}/Z(D)\leq (U\cap 
D)Z(D)/Z(D),$ so  $U_{i}\leq (U\cap 
D)Z(D).$ But then $U_{i}'\leq  ((U\cap 
D)Z(D))'\leq U\cap D.$  By hypothesis, $UU_{i}'/U_{i}'=U/U_{i}'$ is $\sigma$-quasinormal  in  
$G/U_{i}'$ and so $U$ is  $\sigma$-quasinormal in $G$ by Lemma 2.1(3), a
 contradiction.
 Therefore Statement (1) holds.

(2) Assume that this is false.  Let $N=
U^{{\mathfrak{N}}}$ be the nilpotent residual of $U$.
 Then $N < U$ since $U$ supersoluble, so $N$ is $\sigma$-quasinormal in $G$ by 
the minimality of $U$.  
  It is also clear that 
  every proper  subgroup $S$ of $U$ with 
$N\leq S$  is $\sigma$-subquasinormal in $G$, so  $S$ is $\sigma$-quasinormal in $G$.
 Therefore, if $U$ has at least two distinct
  maximal subgroups $S$ and $W$
 such that $N\leq S\cap W$, then $U=\langle S, W \rangle $ is $\sigma$-quasinormal in 
$G$  by Lemma 2.1(4),
 contrary to the  choice of  $U$.
 Hence $U/N$ 
is a cyclic $p$-group for some prime $p$ and $N\ne 1$ since $U$ is not cyclic.

Now we show that $U$ is  a $PT$-group. Let $S$ be a proper 
subnormal subgroup of $U$.  Then $S$ is   $\sigma$-subquasinormal in $G$, so
  $S$ is $\sigma$-quasinormal in $G$ and hence $S$ is $\sigma$-quasinormal in $U$ by 
Lemma 2.1(1).  Therefore $S$ is 
quasinormal in $U$ by Theorem~A. Therefore 
 $U$ is a soluble $PT$-group, so $N=U^{{\mathfrak{N}}}$ is a 
Hall abelian   subgroup of $U$ by\cite[Theorem 2.1.11]{prod}.

It follows that $N\leq U^{{\frak{A}^{*}}}$ and so  $U^{{\frak{A}^{*}}}=NV,$ where 
 $V$ is a maximal
 subgroup of a cyclic Sylow $p$-subgroup $P\simeq U/N$ of $U$.  
Hence  
$NV$ is $\sigma$-quasinormal in $G$ and $NV$ is subnormal in $G$ by  Lemma 2.3(4). 
 Therefore 
$NV$ is quasinormal in $G$ by  Theorem~A.  Assume that for some minimal normal
 subgroup $R$ of $G$ we have  $R\leq (NV)_{G}$. Then $U/R$ is $\sigma$-quasinormal
 in $G/R$
 by hypothesis, so  $U$ is $\sigma$-quasinormal in $G$, a contradiction.  
 Therefore 
$(NV)_{G}=1$, so $NV$ is nilpotent  and  $NV\leq Z_{\infty}(G)$ 
 by \cite[Corollary 1.5.6]{prod} and   then $U=NP$ is nilpotent, so $N=1$, a contradcition.
  Therefore Statement  (2) holds.   
The lemma is proved.

{\bf Lemma 2.7.}  {\sl Suppose that a soluble group $G=D\rtimes M$ satisfies 
  Conditions (i), (ii) and (iii) in Theorem C. If  $A$ 
is a $\sigma$-primary $\sigma$-subnormal subgroup of $G$ such that  
  $A \leq M$, then $D\leq C_{G}(A)$. }

{\bf Proof.}  Let $A$ be a $\sigma _{i}$-group  and $x$  an element 
  of prime power  order $p^{n}$ of $D$.  Let $H_{k}$ be a Hall
 $\sigma _{k}$-subgroup of $G$. Then, by hypothesis, 
$H_{k}=O_{\sigma _{k}}(D)\times S_{k}$, where $  O_{\sigma _{k}}(D)$ and $S_{k}$ are 
 Hall  subgroups of $G$.

  Since  $A$ is 
 $\sigma$-subnormal in $G$, $A \leq H_{i}$ by Lemma 2.6(7) in \cite{1}. On the other
 hand, since 
$A \leq M$, $A\cap D=1$. Therefore 
 $A=(A\cap O_{\sigma _{i}}(D)) \times (A\cap S_{i})=A\cap S_{i}$, 
 so  $A\leq S_{i}$ and hence $ O_{\sigma _{i}}(D)\leq C_{G}(A)$.

 Now, let $k\ne i$.
 Then
$A$ is a Hall $\sigma _{i}$-subgroup of $V:=O_{\sigma _{k}}(D)A$ and $A$ is
 $\sigma$-subnormal
 in $V$ by Lemma 2.6(1) in \cite{1}, so $V=O_{\sigma _{k}}(D)\times A$  by 
Lemma 2.6(10) in \cite{1}   and hence   $ D\leq C_{G}(A)$.  
The lemma is proved.

{\bf Lemma 2.8}   (See Lemma 5.1.9 in \cite{Schm}).
 {\sl Let $A$ be a subgroup of prime power order of $G$. }

(1)  {\sl If $A$ is modular but not subnormal in $G$, then 
$$G/A_{G}=A^{G}/A_{G}\times K/A_{G},$$ where $A^{G}/A_{G}$ is a non-abelian
 $P$-group of order prime to $|K/A_{G}|$.   }

(2)                {\sl   $A$  is modular 
in $G$ if and only if $A$ is modular in $\langle x, A\rangle$ for all  $x\in G$ 
 of prime power order.  }

{\bf Lemma 2.9.} {\sl If $G/Z$ is $p$-closed for some prime $p$ and
 $Z\leq Z_{\infty}(G)$, then $G$ is  $p$-closed.}

{\bf Proof. }  Since $Z\leq Z_{\infty}(G)$, for a
 Sylow $p$-subgroup $Z_{p}$ of $Z$ we have 
   $Z= Z_{p}\times W$, where $Z_{p}$ and $W$
 are characteristic in $Z$ and so normal in $G$.

 Let $P/Z$ be a normal Sylow $p$-subgroup of $G/Z$ and $V$ a Sylow 
 $p$-subgroup of $P$. Then  $Z_{p}\leq V$    and $P=VZ= V\times W$  since 
$W\leq Z_{\infty}(G)\cap P  \leq Z_{\infty}(P)$. Therefore $V$ is characteristic in $P$ and
 so normal in $G$. The lemma is proved.

{\bf Lemma 2.10.}  {\sl Let $G= Q\rtimes P$ be
 a non-abelian $P$-group of type $(q, p)$. }

(1) {\sl $P^{G}=G$.}

(2) {\sl $G/N$ is a non-abelian $P$-group of type $(q, p)$ for
 every proper normal subgroup $N$ of $G$.}

{\bf Proof.}   See Lemma 2.2.2 in \cite{Schm}.

{\bf Lemma 2.11.} {\sl If $A$ and $B$ are normal subgroups of $G$, then every chief
 factor $H/K$  of $G$ below  $AB$ is $G$-isomorphic to either a chiew factor of
 $G$ below $A$ or a chief factor of $G$ between $B\cap A$ and $B$.
}

{\bf Proof.}  This assertion follows from the $G$-isomorphism $AB/A\simeq B/(B\cap A)$ and  
 the Jordan-H\"{o}lder theorem for  the $\Omega$-composition seties of
 a group (see \cite[Chapter A, 3.2]{DH}).

From  Proposition 2.2.8 in \cite{prod}  we get the following useful lemma.

{\bf Lemma 2.12.}  {\sl Let  $\frak{F}$ be a non-empty hereditary formation.}

 (1) {\sl If 
$N$ is a normal subgroup of $G$, then
 $(G/N)^{\frak{F}}=G^{\frak{F}}N/N.$
  }

(2) {\sl If 
$E$ is a  subgroup of $G$, then
 $E^{\frak{F}}\leq G^{\frak{F}}$ and   $N(NE)^{\frak{F}}=NE^{\frak{F}}$. }

{\bf Lemma 2.13.} {\sl  Let $(D, Z(D); U_{1},  \ldots , U_{k})$  be
 a Robinson complex of $G$
and $N$ a normal subgroup of $G$.   }

(1) {\sl $U_{i}'/(U_{i}'\cap Z(D))$ is a simple non-abelian group and 
 $U_{i}'\cap Z(D)=\Phi (U_{i}')=Z(U_{i}')$. In particular, $U_{i}'$ is a
 quasi-simple group. }

(2) {\sl If $N=U_{i}'$ and $k\ne 1$, then $(D/N, Z(D/N);
 U_{1}N/N, \ldots  , U_{i-1}N/N, U_{i+1}N/N,  \cdots , NU_{k}/N)$ is
 a Robinson complex of $G/N$ and   $U_{i}/N=Z(D)N/N=Z(D/N)$.}

(3) {\sl If $N$ is  nilpotent, then $(DN/N, Z(D)N/N;
 U_{1}N/N, \ldots  , NU_{k}/N)$ is
 a Robinson complex of $G/N$ and   $Z(D)N/N=Z(D/N)$.}

(4) {\sl If $p\in \pi (Z(D))$, then $p\in \pi (Z(U_{i}'))$ for some $i$. In particular, 
$p\in \{2, 3\}$.}

{\bf Proof.} Let $Z:=Z(D)=\Phi (D)$.  (1)   
First observe that $U_{i}=U_{i}'Z=U_{i}^{\mathfrak{S}}Z$,
 where  $\mathfrak{S}$ is the class of all soluble groups, 
 since $U_{i}/Z$ is a simple non-abelian group and  so  $U_{i}^{\mathfrak{S}}\leq 
U_{i}' \leq U_{i}^{\mathfrak{S}}$. Hence $U_{i}^{\mathfrak{S}}=U_{i}'$ is perfect.
 On the other hand, 
$U_{i}/Z=U_{i}'Z/Z\simeq U_{i}'/(U_{i}'\cap Z)$  is a simple non-abelian group.
Therefore  $U_{i}'\cap Z=\Phi (U_{i}')=Z(U_{i}')$ since $\Phi (U_{i}')\leq 
\Phi (D)$.

(2), (3) See Remark 1.6.8 in \cite{prod} or Lemma 3.1  in \cite{????}.

(4) Assume that  $p\not \in \pi (Z(U_{i}'))$ for all $i$ and let $Z=Z_{p}\times V$,
 where 
$Z_{p}$ is the Sylow $p$-subgroup of $Z$. Then $Z_{p}\cap U_{i}'=1$, so
$U_{i}'\cap Z=U_{i}'\cap V= \Phi (U_{i}')=Z(U_{i}')$  for all $i$.
 On the other hand, $D=U_{1}  \cdots  U_{k}=ZU_{1}'
  \cdots  U_{k}'=Z_{p}(V(U_{1}'  \cdots  U_{k}')$, so $D=V(U_{1}'  \cdots 
 U_{k}')$ since 
$Z_{p}\leq \Phi (D)$. Hence $Z\leq V(U_{1}'  \cdots  U_{k}')$. But $V$ and 
every subgroup  $ U_{i}'$  has no  a composition factor of order $p$ by Lemma 2.11,
 a contradiction.  Therefore $p\in \pi (Z(U_{i}'))$ for some $i$, where  
 $U_{i}'$ is a quasi-simple group by Part (1). But then $|Z(U_{i}')|$
 divides the order of the Schur multiplier $M(U_{i}'/Z(U_{i}'))$ of $U_{i}'/Z(U_{i}')$.
 Hence $\pi (Z(U_{i}')) \subseteq  \{2, 3\}$ (see Section 4.15(A)
  in   \cite[Ch. 4]{GorI}. Therefore  $p\in \{2, 3\}$.
  Hence we have (4).   
The lemma is proved.

{\bf Lemma 2.14.} {\sl Let $U$ and $N \trianglelefteq  G$ be subgroups of $G$, where
 $U$ is of prime power order. Suppose that $UN/N$ is a modular
 non-subnormal subgroup of  $G/N$.   
Then $$G/(UN)_{G}\simeq U^{G}N/(UN)_{G}\times K/(UN)_{G},$$ where $U^{G}N/(UN)_{G}$
 is a non-abelian $P$-group  of order prime to $|K/UN_{G}|$.  }

{\bf Proof.}  The subgroup $UN/N\simeq U/(U\cap N)$ of $G/N$  is of prime power order,
 so   we can apply Lemma 2.8(1).

  First observe that  $(UN/N)_{G/N}=(UN)_{G}/N$ and
 $(UN/N)^{G/N}=(UN)^{G}/N=U^{G}N/N$. Therefore, by Lemma 2.8(1),   
$$G/(UN)_{G}\simeq    (G/N)/((UN)_{G}/N)=(G/N)/(UN/N)_{G/N}$$$$=   
(UN/N)^{G/N}/(UN/N)_{G/N} \times (K/N)/(UN/N)_{G/N}$$$$=
(U^{G}N/N)/((UN)_{G}/N)  \times (K/N)/((UN)_{G}/N)\simeq
 U^{G}N/(UN)_{G}\times K/(UN)_{G},$$  where 
$$(UN/N)^{G/N}/(UN/N)_{G/N}\simeq U^{G}N/(UN)_{G}$$  is a non-abelian
 $P$-group of order prime to $|(K/N)/(UN/N)_{G/N}|$ and so to
 $|K/(UN)_{G}|$.   
The lemma is proved.

A group $G$ is called \emph{$\pi$-closed} if $G$ has a normal Hall
$\pi$-subgroup.

The following lemma is well-known \cite[Chapter A, 13.2]{DH}.

{\bf Lemma 2.15. } {\sl  If  $H$ is a normal
subgroup of $G$ and  $H/(H\cap \Phi (G))$ is $\pi$-closed, then
$H$ is $\pi$-closed.    }

Recall that a group $G$ is said to be a \emph{$P^{*}$-group} if  $G=A\rtimes \langle t \rangle$,
 where $A$ is an  elementary 
abelian subgroup of $G$, $|t|=r^{n}$ for some prime $r$ and $t$ induces a power
 automorphism of prime order on $A$  \cite[p. 69]{Schm}.

{\bf Lemma 2.16. } {\sl  Let $G=A\rtimes \langle t \rangle$ be a  $P^{*}$-group and let 
$|\langle t \rangle|=p^{n}.$  Then $Z(G)=\langle t^{p} \rangle=\Phi (G)$, $G/Z(G)$ is 
a non-abelian $P$-group and the lattice ${\cal L}(G)$ is modular.}

The following lemma is a corollary of Theorem C.

{\bf Lemma 2.17. } {\sl  If $G$ is a soluible $Q\sigma T$-group, then
 every subgroup of $G$ is a  $Q\sigma T$-group.}

 \section{Proofs of Theorems E and F}

{\bf Proof of Theorem E.}
 First assume that $G$ is a $Q\sigma T$-group.
 Then $G$ is a $PT$-group and every quotient $G/N$ is 
 a $Q\sigma T$-group by Lemma  2.4.  Let
 $D$ be the soluble residual of $G$. Then $D$ is perfect and 
$G/D$ is a soluble group $Q\sigma T$-group, so
 Condition (i) holds for $G$.

Now assume that $D\ne 1$.  Then, in view of  Theorem D,  
$G$   has a Robinson complex
 $(D, Z(D); $ $U_{1},  \ldots , U_{k})$ such that     
 for any set  $\{i_{1}, \ldots , i_{r}\}\subseteq \{1, \ldots , k\}$, where
 $1\leq r  < k$,  $G$ and $G /U_{i_{1}}'\cdots U_{i_{r}}'$ satisfy
 ${\bf N}_{p}$ for all $p\in \{2, 3\} \cap \pi (Z(D))$ and
 ${\bf P}_{p}$  for all $p \in \pi (D)$.  Moreover, in view of  Lemma 2.5, 
 $G$
 and $G /U_{i_{1}}'\cdots U_{i_{r}}'$ satisfy ${\bf Q }_{\sigma P}$ and, in particular, 
satisfy ${\bf Q }_{\sigma (p, q)}$
  for all pairs
 $\{p, q\}\cap  \pi (D)\ne \emptyset$. 
 Hence Conditions (ii) and (iii) hold for $G$.   
  Therefore the necessity
 of the condition of the theorem holds.

 Now, assume, arguing 
by contradiction, that $G$ is a non-$Q\sigma T$-group of minimal order 
satisfying Conditions (i),  (ii), and (iii).

 Then, in view  of Lemma 2.13(4) and  Theorem D,
 $G$ is a $PT$-group, so 
 $D\ne 1$  
and   $G$ has a $\sigma$-subquasinormal $U$ such that $U$ is not $\sigma$-quasinormal
 in $G$ 
but  every $\sigma$-subquasinormal subgroup $U_{0}$   of $G$   with $U_{0} < U$  
 is $\sigma$-quasinormal
 in
$G$.   Let $Z=Z(D).$

(1) {\sl   If   $N$
 is either a non-identity normal nilpotent subgroup of $G$ or   $N=U_{i}'$ for some $i$,
 then $G/N$ is a $Q\sigma T$-group.}

First assume that $k=1$ and $N=U_{1}'$. Then $D'=D=U_{1}=U_{1}'=N$. Therefore 
$G/N=G/D$ is a $Q\sigma T$-group by Condition (i).

Now assume that $k > 1$ and  $N=N_{i}'$.   We 
 can assume  without loss of generality that $i=1$.
    Then
  $ (G/N)/(D/N)\simeq G/D$ is a aoluble $Q\sigma T$-group   and $
(D/N)'=D '/N=D/N$. From Lemma 2.13(2)  it follows that   
 $(D/N, Z(D/N); U_{2}N/N, \ldots , U_{k}N/N)$ is
 a Robinson complex of $G/N$ and 
  $U_{1}/N=ZN/N=Z(D/N)$, where $ZN/N\simeq Z/(Z\cap N)$.  
   Moreover,  by Condition (iii), if
 $\{i_{1}, \ldots , i_{r}\}\subseteq \{2, \ldots , k\}$, where $2\leq r  < k$,
 then the quotients
 $G/N=G/U_{1}'$ and, in view o lemma 2.12(2), $$(G/N) /(U_{i_{1}}N/N)'\cdots (U_{i_{r}}N/N)'=
(G/N)/(U_{i_{1}}'\cdots U_{i_{r}}'U_{1}'/N)\simeq G/U_{i_{1}}'\cdots U_{i_{r}}'U_{1}'$$
 satisfy
 ${\bf N}_{p}$ for all $p\in \{2, 3\} \cap \pi (ZN/N)$,    
 ${\bf P}_{p}$ for all
 $p\in  \pi (D/N))\subseteq \pi (D)$, and  ${\bf Q }_{\sigma (p, q)}$
  for all pairs
 $\{p, q\}\cap  \pi (D/N)\ne \emptyset$.    
 Therefore the hypothesis holds for $G/N=G/U_{1}'$, so $G/U_{i}'$
 is a $Q\sigma T$-group for all $i$ by the choice of $G$.

Similary, it can be proved that if  $N$
 is a non-identity nilpotent normal  subgroup of $G$, then 
 the hypothesis holds for $G/N $ and so $G/N$  is a $Q\sigma T$-group.

(2) {\sl   $U$ is supersoluble.   }

 It is clear that $D=G^{\mathfrak{S}}=G^{\mathfrak{U}}$ and 
$UU_{i}'/U_{i}'$ is $\sigma$-subquasinormal  in  $G/U_{i}'$ by Lemma 2.3(2).
Therefore  $UU_{i}'/U_{i}'$   is  $\sigma$-quasinormal  in $G/U_{i}'$ by Claim (1)
 for all $i$.   Hence $U$ is  supersoluble by Lemma 2.6(1).

(3) {\sl   Suppose that  $N$ is either  
a non-identity normal nilpotent subgroup of $G$ or  $N=U_{i}'$ for
 some $i$.
If    $X$ is a subgroup of $G$ such that $ XN/N$  is $\sigma$-subquasinormal
in $G/N$, then
 $XN/N$ is $\sigma$-quasinormal
 in $G/N$ and $XN$ is $\sigma$-quasinormal in $G$.
In particular,  $U_{G}=1$.}

In view of Claim (1), $G/N$ is a $Q\sigma T$-group, so  $XN/N$ is $\sigma$-quasinormal
 in $G/N$ and hence $XN$ is $\sigma$-quasinormal in $G$ by Lemma 2.1(3). Therefore,
 since $U$ is supersoluble by Claim (2), the choice of $U$ implies that $U_{G}=1$.

  (4) {\sl $U$ is a cyclic $p$-group for some prime $p$ and 
 $V:=U\cap  Z_{\infty}(G)$ is the maximal subgroup of $U$.}

Let $N$  be  a nilpotent 
 non-identity  normal subgroup of $G$. Then $UN/N$ is $\sigma$-subquasinormal
 in $G/N$   by Lemma 2.3(2), so $UN/N$ is $\sigma$-quasinormal
 in $G/N$ by Claim (1). Hence $U$ is a cyclic
 $p$-group for some prime $p$ by Claim (2) and Lemma 2.6(2).

Now,
 let $V$ be the maximal subgroup of $U$. Then  $V=U^{{\frak{A}^{*}}}$  
  is subnormal in $G$ by Lemma 2.3(4) since $U$ is a cyclic
 $p$-group, hence  $V$ is quasinormal in $G$ since $G$ is a $PT$-group. 
  Therefore 
 $V\leq Z_{\infty}(G)$ by \cite[Corollary 1.5.6]{prod} since $V_{G}=1=U_{G}$
 by Claim (3).

(5) {\sl $G$ has a normal subgroup $C_{q}$ of order  $q$ for some  prime $q$.}

If  $Z\ne 1$, it is clear. Now assume
 that $Z=1$. Then  $D= U_{1}\times \cdots 
\times U_{k}$, where $U_{i}$ is a simple non-abelian minimal normal subgroup of $G$
 for all $i$.

Let $E=U_{i}U$, where  $U_{i}\nleq U $.  We show that $U_{i}\leq N_{G}(U) $.   
 Since  $U$ is a  $\sigma$-subquasinormal
subgroup of   $E$  by Lemma 2.3(1), there is  a subgroup chain 
$U=E_{0} < E_{1} < \cdots < E_{t-1} < E_{t}=E$ such that $E_{i-1}$  is a maximal
$\sigma$-quasinormal  subgroup of  $E_{i}$ for all $i=1, \ldots, t$
 and for $M=E_{t-1}$ we have  $M=U(M\cap U_{i})$.  Moreover,  $M\cap U_{i}$ is 
$\sigma$-subquasinormal in $E$ 
  and  $M\cap U_{i} < U_{i}$ since 
$M < E$, so $M\cap U_{i}=1$. Therefore $U=M$ is  a maximal
$\sigma$-quasinormal  subgroup of  $E$. Assume that $U$ is not normal in $E$. Then, 
 by  Lemma 2.2, $E/U_{E}$ is a group of order $qr$ for primes  $q$ and $r$. But this is
 imposible since 
 $U_{i}\simeq U_{i}U_{E}/U_{E}\leq E/U_{E}$.  
 Therefore  $U_{i}\leq N_{E}(U)$ for all $i$,
 so 
$D\leq N_{G}(U)$ and hence $U\cap D\leq O_{p}(D)=1$ by Claim (4).

 It follows than $DU=D\times U$, so 
$1  < U\leq C_{G}(D)$. But  $C_{G}(D)\cap D=1$ since $Z=Z(D)=1$.
 Therefore $C_{G}(D)\simeq C_{G}(D)D/D$ is soluble. Hence for some prime $q$
 dividing $|C_{G}(D)|$ we have 
$O_{q}(C_{G}(D))\ne 1$. But $O_{q}(C_{G}(D))$ is characteristic in $C_{G}(D)$,
so 
$O_{q}(C_{G}(D))$ is normal in $G$ and hence  we have (5).

(6) {\sl $U^{G}$ is soluble}.

The  subgroup $C_{q}U/C_{q}$  is $\sigma$-subquasinormal
 in $G/C_{q}$ by Lemma 2.3(2), so this subgroup is $\sigma$-quasinormal  and 
hence  modular in  $G/C_{q}$ by Claim (3).

 First assume that $C_{q}U/C_{q}$  is not subnormal  in $G/C_{q}$.
 Then, by Lemma 2.14,  
$C_{q}U^{G}/(C_{q}U)_{G}$ 
  is a non-abelian   $P$-group, so $C_{q}U^{G}/(C_{q}U)_{G}$ 
 is soluble and hence  
  $$U^{G}(C_{q}U)_{G}/(C_{q}U)_{G}\simeq U^{G}/(U^{G}\cap   (C_{q}U)_{G})$$ is soluble since
 $C_{q}U$ is soluble.
 Hence $U^{G}$  is soluble.

Now assume that  $C_{q}U/C_{q}$  is  subnormal  in $G/C_{q}$,  so
 $$U^{G}/(U^{G}\cap C_{q})\simeq
 C_{q}U^{G}/C_{q}=(C_{q}U/C_{q})^{G/C_{q}}\leq O_{p}(G/C_{q})$$  by Claim (4).
 Hence  $U^{G}$ is soluble.

(7) {\sl $U$ is not subnormal in $G$. }

  Assume that   $U$ is subnormal in $G$. 
 Then $U$ is quasinormal and
 so $\sigma$-quasinormal  in $G$ since 
$G$ is a $PT$-group,  a contradiction. Hence we have (7).

(8)  $|U|=p$.

 Assume that  $|U| > p$. Then $1 < V\leq R:=O_{p}(Z_{\infty}(G))$ by Claim (4) and 
  $U\nleq R$ by Claim (7). Denote $E=RU$.   Then $E^{G}=U^{G}R$ and, in view
 of Claims (4) and (7),  $E$  is not subnormal in $G$. Moreover, $E$ 
 is $\sigma$-quasinormal and so
 modular in $G$ by Claim (3).  The group 
  $UR/R\simeq U/(U\cap R)=U/V$  has order $p$, so $(RU)_{G}=R$.
   
 In view of Lemma 2.14,
 $$G /R=G/E_{G}\simeq  E^{G}/E_{G}\times K/E_{G}=U^{G}R/R\times K/R,$$
 where 
 $RU^{G}/R\simeq U^{G}/(U^{G}\cap R)$ 
 is a   non-abelian  $P$-group of order prime to $K/R$.

The group  $ U^{G}/(U^{G}\cap R)$ is  $q$-closed for some prime $q$
 dividing its order and 
so $ U^{G}$ is  $q$-closed by Lemma 2.9 since $U^{G}\cap R\leq  Z_{\infty}(U^{G})$. 
If $Q$ is the  normal Sylow $q$-subgroup of $ U^{G}$, then $U\nleq Q$ by Claim (7), so 
$q\ne p$. Therefore $ U^{G}/(U^{G}\cap R)$ 
  is a   non-abelian  $P$-group of type $(q, p)$, so  $U^{G}=Q\rtimes P$, where 
$P$ is a non-normal Sylow  $p$-subgroup of $U^{G}$, containing $U$, and $Q$ is
 subnormal in $G$.    
In particular, $U^{G}$ and  $RU^{G}/R$ are $\pi$-groups,
where $\pi =\{p, q\}$, so $G$ is $\pi$-soluble and hence 
  $D$ and $D/Z$ are $\pi$-soluble groups.

 Assume that  $U^{G}\cap D\ne 1$.    Since 
 $U^{G}\cap D\leq Z\leq \Phi (D)$ by Claim (6),
 for some $i$ and for some  $r\in \{p, q\}$ the mumber $r$ divides $|U_{i}/Z|$.
 It follows that $U_{i}/Z$ is an abelian group, a contradiction.

Therefore $U^{G}\cap
 D= 1$  and 
 so     $$U^{G}\simeq U^{G}/(U^{G}\cap D)\simeq U^{G}D/D=
(QD/D)\rtimes  (PD/D),$$  where   
  $G/D$ is a soluble  $Q\sigma T$-group by Condition (i).  
Therefore, in view of Theorem C,   
$G/D=T\rtimes L$, where $T=(G/D)^{\frak{N_{\sigma}}}$  and the following  hold:  
  $T$   is an abelian  Hall
 subgroup of $G/D$ and all subgroups of $T$ are normal $G/D$ and  the  
 lattice ${\cal L}(L)$,
of all subgroups  of $L$,  is modular. Then  $PD/D\nleq T$,                            
 so $UD/D\leq PD/D\leq L^{x}$   for some $x\in G/D$.

First assume that  $QD/D\leq T$. Since $UD/D$ is a
 $\sigma$-subquasinormal $p$-subgroup  of $ G/D$ by Lemmma 2.3(2), 
 $T\leq C_{G/D}(UD/D)$ by Lemma 2.7, so
 $$(QD/D)\rtimes  (PD/D)=U^{G}D/D=(UD)^{G}/D= 
 (UD/D)^{G/D}=(UD/D)^{TL^{x}}=(UD/D)^{L^{x}}\leq L^{x},$$ a contradiction. Hence  		
 $QD/D\nleq T$ and  so  $(QD/D)\rtimes  (PD/D)\leq L^{x}$ since $T$ and $L^{x}$ 
are  Hall subgroups of $G/D$ and $QD/D$ is a  subnormal $q$-subgroup of $G/D$.

In view of   Theorem 2.4.4 in \cite{Schm}, $L^{x}$ is a direct product of
 $P^{*}$-groups $P_{i}$ and primary groups  $Q_{j}$
 (that is, $Q_{j}$ is of prime power order) with relatively prime
 orders.  Then  for any   factor 
$Q_{j}$ of $L^{x} $  we have
 $Q_{j}\leq Z_{\infty}(L^{x})$, so 
 $QD/D\nleq Q_{j}$ and  $PD/D\nleq Q_{j}$  for all $j$ since
 $U^{G}\simeq U^{G}D/D\simeq QD/D\rtimes PD/D$ is not nilpotent.

 Therefore   for some $i$  and $k$ we have 
 $QD/D\leq P_{i}$ and $PD/D\leq P_{k}$, where $[P_{i}, P_{k}]=1$ for $i\ne k$, so $i=k$.
Hence $QD/D\rtimes PD/D\leq P_{i}=A\rtimes \langle t \rangle$, where $A$ is elementary 
abelian subgroup of $P_{i}$, $|t|=r^{n}$ for some prime $r$ and $t$ induces a power
 automorphism of prime order on $A$. Therefore 
$A$ is a $q$-group and $t$ is a $p$-element of $P_{k}$ by Lemma 2.16. Hence $P\simeq 
PD/D$  is a cyclic $p$-group.

 Since $U^{G}/(U^{G}\cap R)$ 
 is a   non-abelian  $P$-group and $U\nleq U^{G}\cap R$,   
$ U(U^{G}\cap R)/(U^{G}\cap R)$ is a Sylow $p$-subgroup of 
$U^{G}/(U^{G}\cap R)$ and  so  
  $U(U^{G}\cap R)$ is a cyclic Sylow $p$-subgroup of $U^{G}$.
It follows that either  $U(U^{G}\cap R)=U$ or $U(U^{G}\cap R)=U^{G}\cap R$.
 In the former case
 we have $U^{G}\cap R=V\leq U_{G}$, which is impossible by  Claim (3), 
 so 
$U(U^{G}\cap R)=U^{G}\cap R$ and hence $U$ is subormal in $G$, contrary to Claim (7).
Therefore  we have (8).

(9) {\sl $U\nleq D$.}

Assume $U\leq D$. 
From Claim (7) it follows that $U\nleq Z$ and then, by Claim (8) and
 Lemma 2.3(1)(2)(5), for some $i$ we have 
 $U\simeq UZ/Z=U_{i}/Z$,  a
 contradiction.  Hence we have~(9).

(10) $O_{p}(D)=1$.

Assume that $G$ has a normal subgroup $Z_{p}\leq Z=\Phi(D)$ of order $p$.
Then $Z_{p}U$ is not subnormal in $G$ by Claim (7) and, also,  $(Z_{p}U)_{G}=Z_{p}$
 by Claim (8) and 
$(Z_{p}U)^{G}=Z_{p}U^{G}$, so 
$G/Z_{p}\simeq Z_{p}U^{G}/Z_{p}\times K/Z_{p},$ where $Z_{p}U^{G}/Z_{p}$ is
 a non-abelian $P$-group of order $p^{a}q^{b}$ prime to $|K/Z_{p}|$ by Lemma 2.14. Hence 
$G/Z_{p}$,   $D/Z_{p}$, and $D$  are $\{p, q\}$-soluble and  $p$ divides $|D/Z_{p}|$. 
Hence $O_{p}(D/Z)\ne 1$. This contradiction completes the proof of the  claim.

(11) {\sl $U_{i}'U=U_{i}'\times U$ and so  $U_{i}'U$ is not subnormal in $G$
 for all $i$.  }

In view of Claims (8) and (9), it is enough to show that $U_{i}'\leq N_{G}(U)$. 
 
By Lemma 2.13(1),   $U_{i}'\cap Z=\Phi (U_{i}')=Z(U_{i}')$
  and $U_{i}'/\Phi (U_{i}')$ is a simple non-abelian group. In particular, 
$U_{i}'$  is quasi-simple.

Let $E=U_{i}'U=U_{i}'\rtimes U$. Then $E'=U_{i}'$.   Let 
$U=E_{0} < E_{1} < \cdots < E_{t-1} < E_{t}=E$ be  a   subgroup chain such that
 $E_{i-1}$
  is a maximal $\sigma$-quasinormal  subgroup 
 of  $E_{i}$ for all $i=1, \ldots, t$ and for $M=E_{t-1}$ we have 
$M=U(M\cap U_{i}')$. Then, by Lemma 2.2, we have   either
 $M=E_{t-1}$ is  a maximal normal subgroup of $E$ or $M$
 is a maximal subgroup of $E$ such  that 
$E/M_{E}$ is a $\sigma$-primary non-abelian group of order $qr$
 for some primes  $q$ and $r$.

First assume that $M$ is normal in $E$. 
 From $E=U_{i}'U=U_{i}'M$ it follows 
that   
 $E/M\simeq U_{i}'/(M\cap U_{i}')$ is a
 simple  group and  so $U_{i}'/(M\cap U_{i}')$    is  a
 simple  non-abelian group since  $U_{i}'$ is perfect. Therefore 
 $M\cap U_{i}'=U_{i}'\cap Z =\Phi (U_{i}')$ is a $p'$-group
 by Claim (10), so
 $U$ is a Sylow $p$-subgroup
 of $M=U(M\cap U_{i}')$.
 Then, by the Frattini argument,
 $E=MN_{E}(U)=  (M\cap U_{i}')N_{E}(U)=\Phi (U_{i}')N_{E}(U)$.   
 But  $\Phi (U_{i}')\leq \Phi (E)$, therefore  $N_{E}(U)=E$ and    so $U_{i}'\leq N_{G}(U)$. 
 
Finally, assume that $E/M_{E}$ is a non-abelian group of order $qr$
 with $V/M_{E}=(E/M_{E})'$.
 Then $|(E/M_{E})/(E/M_{E})'|= (E/M_{E})/(V/M_{E})= |E/V|$ is a prime,
so $V=U_{i}'$.   Hence $M_{E}\leq U_{i}'$ and $U_{i}'/M_{E}$ is a non-identity soluble 
group, so $U_{i}'$ is not perfect.  This contradiction shows that  
 we have (11).

(12) {\sl $U^{G}$ is not a non-abelian $P$-group. }

Assume that   $U^{G}$ is a non-abelian $P$-group. Then, in view of Claim (7),
$U^{G}=Q\rtimes U$ is   of type $(q, p)$ for some prime $q$. Let  $\pi =
\{q, p\}$.

 First suppose that  $\pi \cap \pi (D)= \emptyset$.
 Then $U^{G}\cap D=1$, so  $[U^{G}, D]=1$ and $G$ is $\pi$-soluble.

 We show that $G/D$ is  $\pi$-decomposable.  Let 
  $N=U_{1}'$ and $F= NU=N\times U$. Then 
 
 $F^{G}=NU^{G}$ and
 $F_{G}=N$ by Claim (8).   In view of Claims (3) and   (7), 
$F/N$ is not subnormal  but modular in $G/N$ and  so    
$$G/N\simeq NU^{G}/N \times K/N,$$ where 
$NU^{G}/N=O_{\pi}(G/N)$ and $K/N=O_{\pi'}(G/N)$, by Lemma 2.14. Therefore  
 $G/N$ is $\pi$-decomposable. Hence 
$G/D$ is $\pi$-decomposable.

Let $E$ be a minimal supplement to $D$ in $G$. Then $E\cap D\leq \Phi (E)$, so $E$ is
 soluble and $\pi$-decomposable, that is, $E=O_{\pi}(E)\times O_{\pi'}(E)$ by Lemma 2.15
 since
 $G/D\simeq E/(E\cap D)$.

Let $x\in G_{r}$, where $G_{r}$  be a Sylow $r$-subgroup of  $G$. 
Assume that $r\not \in  \pi$. Then for some Sylow
 $r$-subgroup $D_{r}$ of $D$ 
and a Sylow $r$-subgroup $E_{r}$ of $E$ and some $y\in G$ we have 
 $G_{r}=D_{r}E_{r}^{y}$.

 Hence $x=de$, where  $d\in D_{r}$ and $e\in E_{r}^{y}$. Then  $d\leq C_{G}(U)$  
since $[U^{G}, D]=1$.  Since $|G:E_{r}^{y}|=|DE_{r}^{y}:E_{r}^{y}|=|D:D\cap E_{r}^{y}|$ is 
a $\pi'$-mumber,  the  Hall $\pi$-subgroup $O_{\pi}(E^{y})$ of $E^{y}$ 
 is a Hall $\pi$-subgroup 
of $G$. Hence $U^{G}\leq O_{\pi}(E^{y})$ and so $e\leq C_{G}(U)$ since
 $e\in O_{\pi'}(E^{y})$. 
 Therefore  $x\leq C_{G}(U)$ and hence $U$ is normal  and so modular
 in 
 $\langle x, U \rangle$.

Now let  $r\in \pi$.   Then  $V=U^{G}G_{r}$ is a $\pi$-subgroup  $G$, so 
  $V\cap D=1$ and therefore $V\simeq VD/D$  is a soluble $Q\sigma T$-group by
 Lemma 2.17,
 so  $U$ is  $\sigma$-quasinormal and so modular in $V$. Hence 
$U$ is modular in $\langle x, U \rangle$ by \cite[Page 201, Property  (2)]{Schm}. 
  Therefore $U$ is modular in $G$ by  Lemma 2.8(2)
 and so $U$ is $\sigma$-quasinormal in $G$,   a contradiction.

 Finally, if 
  $\pi \cap \pi (D)\ne  \emptyset$, then $G$  satisfies   
 ${\bf Q }_{\sigma (p, q)}$ by Condition (iii), so $U$ is modular and so 
 $\sigma$-quasinormal
   in $G$.
This contradiction completes the proof of the claim.

 (13) {\sl  $G$ has a normal subgroup $C_{q}$ of prime
 order $q$ such that $C_{q}\leq Z(U_{1}')=\Phi (U_{1}')$.}

 Let $E=U_{1}' U=U_{1}'\times U$.  Then $E$ is modular and 
 not subnormal in $G$ by  Claims   (3) and (7). Moreover, $E_{G}=U_{1}'$ by Claim (8)
 and
 $E/U_{1}'\simeq U$ is  a modular non-subnormal subgroup of $G/U_{1}'$. Hence 
$$E^{G}/E_{G}=U^{G}U_{1}'/U_{1}'\simeq U^{G}/(U^{G}\cap U_{1}')$$ is a
 non-abelian $P$-group by Lemma 2.8(1).
Hence  $1 < U^{G}\cap U_{1}' \leq Z(U_{1}')$ by Claims  (6) and (12). Hence $G$
 has a normal subgroup $C_{q}$ of prime
 order $q$ such that $C_{q}\leq U_{1}'$. But  $U_{1}'$ is a quasi-simple group by Lemma 2.13(1)
 and so 
 $C_{q}\leq Z(U_{1}')=\Phi (U_{1}')$.

 {\sl Final contradiction.}  From Claims  (7), (9) and (11) it follows that
 $E=C_{q}U=C_{q}\times U$ is not subnormal in $G$ and, in view of Claim (8),
 $E_{G}=C_{q}$.
Hence $G/E_{G}\simeq  E^{G}/E_{G}\times K/E_{G},$ where
 $$E^{G}/E_{G}=C_{q}U^{G}/C_{q}\simeq 
U^{G}/(C_{q}\cap U^{G})$$ is a non-abelian
 $P$-group of order prime to $|K/C_{q}|$ by Lemma 2.8(1). Hence $G$ is a $\pi$-soluble 
group, where $\pi= \pi (U^{G}/(C_{q}\cap U^{G}))$. Then $D/Z$ is 
$\pi$-soluble. But $C_{q}\leq \Phi (U_{1}')\leq  \Phi (D)=Z$ by Claim (1),
 so $q$ divides $|D/Z|$. Hence
$q$ does not divides  $|C_{q}U^{G}/C_{q}|$.

If $C_{q}\cap U^{G}=1$, then $U^{G}\simeq  C_{q}U^{G}/C_{q} $ is a non-abelian
 $P$-group, contrary to Claim (12),  so  $C_{q}\leq  U^{G}$. Then   
 $C_{q}$ is a Sylow $q$-subgroup of $U^{G}$.
 Hence $U^{G}=C_{q}\rtimes (R\rtimes U)$,
 where  $R\rtimes U\simeq  U^{G}/C_{q}$ is a non-abelian $P$-group. 
Let $C=C_{U^{G}}(C_{q})$. Then $U\leq C$ by Claim (11) and so, by Lemma 
 2.10(1),  $R\rtimes
 U=U^{R\rtimes U}\leq C$. Hence $C_{q}\leq Z(U^{G})$.
 Therefore
  $U^{G}=C_{q}\times (R\rtimes U)$, where $R\rtimes U$ is characterisric in $U^{G}$
 and so it is normal in $G$. But then $U^{G}=R\rtimes U\ne C_{q}\rtimes (R\rtimes U)$, 
a contradiction.             
The theorem is proved.  

{\bf Proof of Theorem F.}  In view of Example 1.2(i), 
Teorem F is a special case of Theorem E, where  $\sigma =\{\mathbb{P}\}$.

\section{Final remarks, further applications}

1. First Consider the special case of Theorem E where
 $\sigma =\sigma ^{1\pi }=\{\{p_{1}\}, \ldots , \{p_{n}\}, \pi'\}$ and   
$\pi= \{p_{1}, \ldots , p_{n}\}$ (see Example 1.2(iii)). 

In this case  we say that  
   $G$ is a \emph{$Q 1\pi  T$-group}
 if  $1\pi $-quasinormality    
is a transitive relation on $G$, 
 and we also say in this case  that "$G$  satisfies ${\bf Q }_{1\pi  (p, q)}$"
 instead of "$G$ satisfies ${\bf Q }_{\sigma  (p, q)}$".

Observe that $G$ satisfies ${\bf Q }_{1\pi  (p, q)}$       
if  whenever $N$ is  a soluble normal
subgroup of $G$ and $P/N$ is a normal non-abelian 
 $P$-subgroup  of type  $(p, q)$  of $G/N$,  where $p, q\in \pi' $,
 every   subgroup of $P/N$ is
 modular in   $G/N$.  
Therefore we get from Theorem E the following result.

{\bf Corollary 4.1.}  {\sl A group $G$ is a $Q 1\pi  T$-group if   
 and   only if  $G$  has a perfect normal subgroup $D$ such that:}

(i) {\sl  $G/D$ is a soluble  $Q 1\pi  T$-group,  }

(ii) {\sl if  $D\ne 1$,  $G$ has a Robinson complex
 $(D, Z(D); U_{1},  \ldots , U_{k})$ and }

(iii) {\sl   for any set $\{i_{1}, \ldots , i_{r}\}\subseteq \{1, \ldots , k\}$, where
 $1\leq r  < k$,  the groups $G$ and $G /U_{i_{1}}'\cdots U_{i_{r}}'$ satisfy
 ${\bf N}_{p}$ for all $p\in \{2, 3\}\cap  \pi (Z(D))$,
 ${\bf P}_{p}$ 
for all $p\in \pi (D)$, and ${\bf Q }_{1\pi (p, q})$
  for all pairs  $\{p, q\}\subseteq \pi'$ with 
 $\{p, q\}\cap  \pi (D)\ne \emptyset$.}

2. Now Consider the special case of Theorem E where
 $\sigma =\sigma ^{\pi }=\{\pi, \pi'\}$ (see Example 1.2(iv)). 

In this case  we say that  
   $G$ is a \emph{$Q \pi, \pi'  T$-group}
 if  $\pi, \pi' $-quasinormality     
is a transitive relation on $G$, 
 and we also say in this case that $G$  "satisfies ${\bf Q }_{\pi, \pi'  (p, q)}$"
 instead of "$G$ satisfies ${\bf Q }_{\sigma  (p, q)}$".

Observe that $G$ satisfies ${\bf Q }_{\pi, \pi'  (p, q)}$       
if  whenever $N$ is  a soluble normal
subgroup of $G$ and $P/N$ is a normal non-abelian 
 $P$-subgroup  of type  $(p, q)$  of $G/N$,  where $p, q\in \pi _{0}\in \{\pi, \pi'\} $,
 every  subgroup of $P/N$ is  modular in $G/N$.

Therefore we get from Theorem E the following result.

{\bf Corollary 4.2.}  {\sl A group $G$ is a $Q \pi, \pi'  T$-group if   
 and   only if  $G$  has a perfect normal subgroup $D$ such that:}

(i) {\sl  $G/D$ is a soluble  $Q\pi, \pi'  T$-group,  }

(ii) {\sl if  $D\ne 1$,  $G$ has a Robinson complex
 $(D, Z(D); U_{1},  \ldots , U_{k})$ and }

(iii) {\sl   for any set $\{i_{1}, \ldots , i_{r}\}\subseteq \{1, \ldots , k\}$, where
 $1\leq r  < k$,  the groups $G$ and $G /U_{i_{1}}'\cdots U_{i_{r}}'$ satisfy
 ${\bf N}_{p}$ for all $p\in \{2, 3\}\cap \pi (Z(D))$,
 ${\bf P}_{p}$ 
for all $p\in \pi (D)$, and ${\bf Q }_{\pi, \pi' (p, q})$ for all pairs 
  $\{p, q\}\cap  \pi (D)\ne \emptyset$.}

3.  In the  case  when   $\sigma =
\sigma ^{1}=\{\{2\}, \{3\}, \{5\} \ldots  \}$  (see Example 1.2(ii))
we get from Theorem E the following clarification of Theorem D.

{\bf Corollary 4.3.} 
 {\sl $G$ is a $PT$-group if 
 and   only if  $G$  has a normal perfect subgroup $D$ such that:}

(i) {\sl  $G/D$ is a soluble $PT$-group, and }

(i) {\sl if $D\ne 1$, $G$ has a Robinson complex
 $(D, Z(D); U_{1},  \ldots , U_{k})$ and }

(iii) {\sl   for any set  $\{i_{1}, \ldots , i_{r}\}\subseteq \{1, \ldots , k\}$, where
 $1\leq r  < k$,  $G$ and $G /U_{i_{1}}'\cdots U_{i_{r}}'$ satisfy
 ${\bf N}_{p}$ for all $p\in \{2, 3\}\cap  \pi (Z(D))$ and
 ${\bf P}_{p}$ 
for all $p \in \pi (D)$. } 

  4.   In the paper \cite{Archiv}, the following special case of Theorem F was proved.

{\bf Corollary 4.4.}  {\sl
 A group $G$ is an $MT$-group if   
 and   only if  $G$  has a perfect normal subgroup $D$ such that:}

(i) {\sl  $G/D$ is an $M$-group,  }

(ii) {\sl if  $D\ne 1$,  $G$ has a Robinson complex
 $(D, Z(D); U_{1},  \ldots , U_{k})$ and }

(iii) {\sl   for any set $\{i_{1}, \ldots , i_{r}\}\subseteq \{1, \ldots , k\}$, where
 $1\leq r  < k$,  $G$ and $G /U_{i_{1}}'\cdots U_{i_{r}}'$ satisfy
 ${\bf N}_{p}$ for all $p\in  \pi (Z(D))$,
 ${\bf P}_{p}$ 
for all $p\in \pi (D)$, and ${\bf M}_{p, q}$ for all pairs
 $\{p, q\}\cap \pi (D)\ne \emptyset.$}

{\bf Remark 4.5.}  Theorem F not only strengthens Corollary 4.4 
 but also gives a new proof of it.


\begin{thebibliography}{s2}




\bibitem{Schm} R.~Schmidt, {\it Subgroup Lattices of Groups},
  Walter de Gruyter, Berlin, 1994.



\bibitem{commun}  A.N. Skiba, On some results in the theory of
 finite partially  soluble groups, \emph{Commun. Math. Stat.},
{\bf 4}(3)  (2016), 281--309.


\bibitem{1}  A.N.  Skiba, On $\sigma$-subnormal
 and $\sigma$-permutable subgroups of finite groups,
\emph{J. Algebra}, {\bf 436} (2015), 1--16.


 
                                       



\bibitem{5} O. Ore, Contributions in the theory of groups of finite
order, \emph{ Duke Math.}, {\bf 5} (1939), 431--460.




       


\bibitem{It} N. Ito, J. Sz\'{e}p, Uber die Quasinormalteiler von endlichen
Gruppen, \emph{Act. Sci. Math.}, {\bf 23} (1962), 168--170.




\bibitem{MaierS} R. Maier, P. Schmid, The embedding of permutable
subgroups in finite groups, \emph{Math. Z.}, {\bf 131} (1973), 269--272.



\bibitem{Th}  J.G. Thompson, An example of core-free quasinormal subgroup of $p$-group,
 \emph{Math. Z.}, {\bf 96} (1973), 226--227.


\bibitem{ProblemI} A.N.  Skiba, On $\sigma$-properties   of finite groups I, 
\emph{Problems of Physics, Mathematics and Technics}, 4(21) (2014), 89--96.
        



\bibitem{comm} H. Li, A.-M. Liu, I.N. Safonova, 
 A.N. Skiba Characterizations of some
 classes of finite $\sigma$-soluble  $P\sigma T$-groups, \emph{Communications
 in Algebra},  {\bf  52} (1) (2024), 128--139. https://doi.org/10.1080/ 00927872.2023.2235006.
 
  


\bibitem{???} A-M. Liu,  M. Chen, I.N. Safonova, A.N. Skiba, 
Finite groups with modular  $\sigma$-subnormal subgroups,   
\emph{J. Group Theory.} https://doi.org/10.1515/jgth-2023-0064.
 

 
\bibitem{????}  X.-F. Zhang, W. Guo, I.N. Safonova, A.N.Skiba,  
 A Robinson description of finite $P\sigma T$-groups, 
 \emph{J. Algebra, }  {\bf 631} (2023),    218--235.

                          

\bibitem{KegI} O.H. Kegel,  Untergruppenverbande endlicher Gruppen, die den
 subnormalteilerverband each enthalten, \emph{Arch. Math.}, {\bf 30}(3)
  (1978), 225--228.


\bibitem{Hu11} B. Hu, J. Huang, A.N. Skiba,   On $\sigma$-quasinormal
 subgroups of finite groups,
  \emph{Bull. Austral. Math. Soc.,} {\bf 99}(3) (2019), 413--420.




\bibitem{prod}  A. Ballester-Bolinches, R. Esteban-Romero, M.
  Asaad,  \emph{Products of Finite Groups},  Walter de Gruyter, Berlin-New York,  2010.




\bibitem{zaher}  G. Zacher, I gruppi risolubili finiti in cui i sottogruppi
 di composizione coincidono con
 i sottogruppi quasi-normali. \emph{Atti della Accademia Nazionale dei Lincei Rend. cl. Sci.
 Fis. Mat. Natur.}, (8) 37 (1964), 150--154



\bibitem{217} D.J.S. Robinson, The structure of finite groups in which permutability is a
transitive relation,
 \emph{J. Austral. Math. Soc.}, {\bf 70} (2001), 143--159.





\bibitem{MZ}    X.-F.~Zhang,  W.~Guo,   I.N.~Safonova,  
A.N.~Skiba,  Finite soluble groups with transitive  $\sigma$-quasinormality relation, 
\emph{Math. Notes}, {\bf 114}(5) (2023), 1029--1036.   
https://doi.org/10.1134/ S0001434623110342. 

           \bibitem{A. Frigerio}  A. Frigerio, Gruppi finiti nei quali e transitivo
 l'essere sottogruppi modulare, \emph{ Ist. Veneto Sci. Lett.
	 Arti, Atti Cl. Sci. mat. natur.}, {\bf 132}, (1973--1974), 185--190. 

\bibitem{mod} I. Zimmermann, Submodular subgroups of finite groups, \emph{Math. Z.},
 {\bf 202}  (1989), 545--557.

                  

\bibitem{DH}  K. Doerk, T. Hawkes,  \emph{Finite Soluble Groups},
 Walter de Gruyter, Berlin-New York, 1992.






\bibitem{Archiv} A.-M. Liu, W. Guo, I.N. Safonova, A.N.
 Skiba , Finite groups in which modularity is a transitive relation, 
 \emph{Archiv der Mathematik}, {\bf 121}(2) (2023),  111--121.
                                                   


%\bibitem{Alg2023} Xin-Fang Zhang, Wenbin Guo, Inna N. Safonova,A. Skiba, A Robinson description of finite $P\sigma T$-groups, \emph{Journal of Algebra}, {\bf 631} (2023), 218--235.  


\bibitem{GorI} D. Gorenstein, \emph{Finite simple groups. An introduction to 
their Classification},   Plenum Press New York and London, 1982.

 
                                       




 \end{thebibliography}
\end{document}